\documentclass[aop,seceqn,MSNbibl,citesort,dvips]{arximspdf}
\usepackage{graphicx}

\doi{10.1214/10-AOP568}
\volume{39}
\issue{3}
\pubyear{2011}
\firstpage{881}
\lastpage{903}

\makeatletter

\newtheorem{theorem}{Theorem}[section]
\newtheorem{lemma}[theorem]{Lemma}
\newtheorem{corollary}[theorem]{Corollary}
\newtheorem{proposition}[theorem]{Proposition}
\newproclaim{definition}[theorem]{Definition}

\newcommand{\iint}{\int\!\!\int}
\newcommand{\iiiint}{\int\!\!\int\!\!\int\!\!\int}

\newcommand{\Var}{\operatorname{Var}}
\newcommand{\Cov}{\operatorname{Cov}}
\newcommand{\SL}{\operatorname{SL}}
\newcommand{\Aff}{\operatorname{Aff}}
\newcommand{\Area}{\operatorname{Area}}

\makeatother

\begin{document}
\begin{frontmatter}

\title{Central limit theorems for random polygons in an arbitrary
convex set}
\runtitle{Random polygons}

\begin{aug}
\author[A]{\fnms{John} \snm{Pardon}\corref{}\ead
[label=e1]{jpardon@princeton.edu}}
\runauthor{J. Pardon}
\affiliation{Princeton University}
\address[A]{Department of Mathematics\\
Princeton University\\
Princeton, New Jersey 08544\\
USA\\
\printead{e1}} 
\end{aug}

\received{\smonth{3} \syear{2010}}
\revised{\smonth{5} \syear{2010}}

%
\begin{abstract}
We study the probability distribution of the area and the number of
vertices of random polygons in a convex set $K\subset\mathbb R^2$. The
novel aspect of our approach is that it yields uniform estimates for
all convex sets $K\subset\mathbb R^2$ without imposing any regularity
conditions on the boundary $\partial K$. Our main result is a central
limit theorem for both the area and the number of vertices, settling a
well-known conjecture in the field. We also obtain asymptotic results
relating the growth of the expectation and variance of these two
functionals.
\end{abstract}

%

\end{frontmatter}

\section{Introduction}\label{sec1}

Consider a Poisson point process in a convex set $K\subset\mathbb R^2$
of intensity equal to the Lebesgue measure. We denote by $\Pi_K$ the convex
hull of the points of this process; $\Pi_K$ is called a \textit{random
Poisson polygon}. We denote by $N=N(\Pi_K)$
the number of vertices of $\Pi_K$ and by $A=A(\Pi_K)$ the area of
$K\setminus\Pi_K$. In this paper, we
develop techniques to study the distributions of these random
variables. Our main result is a central limit theorem, which
is \textit{uniform} over the set of all convex $K\subset\mathbb R^2$:
\begin{theorem}\label{strongclt}
As $\Area(K)\to\infty$, we have the following central limit theorems
for $\Pi_K$:
%
%
\begin{eqnarray}
&&\sup_x\biggl|P\biggl(\frac{N-\mathbb E[N]}{\sqrt{\Var N}}\leq
x
\biggr)-\Phi(x)\biggr|\ll\frac{\log^2\mathbb E[N]}{\sqrt{\mathbb
E[N]}},
\\
%
%
&&\sup_x\biggl|P\biggl(\frac{A-\mathbb E[A]}{\sqrt{\Var A}}\leq
x
\biggr)-\Phi(x)\biggr|\ll\frac{\log^2\mathbb E[A]}{\sqrt{\mathbb E[A]}}.
\end{eqnarray}
Here $\Phi(x)=P(Z\leq x)$ where $Z$ is the standard normal distribution.
\end{theorem}

The novel aspect of our approach is that we require no regularity on
$\partial K$; it is this that enables us to obtain bounds which are
uniform over all convex sets. Previous results on random polygons
analogous to Theorems \ref{strongclt}
have been confined to two cases: (i) $K$ a polygon \cite{groen1,groen2}
and (ii) $\partial K$ of class $C^2$ with nonvanishing curvature~\cite{hsing}.
The key part of our argument is our use of a new compactness result for
various types of local configuration spaces of convex boundaries.

As a consequence of our techniques, we also prove the following:
\begin{theorem}\label{strongvarest}
As $\Area(K)\to\infty$, we have the following estimates for $\Pi
_K$\footnote{After this paper was written, we learned that Imre B\'ar\'
any and Matthias Reitzner have independently proved this result, as
well as the closely related Corollary \ref{uniformvarest}.}:
%
%
\begin{equation}
\mathbb E[N]\asymp\Var N\asymp\mathbb E[A]\asymp\Var A.
\end{equation}
\end{theorem}

In other words, there is (up to a constant factor) only one parameter,
say $\mathbb E[A]$,
which controls the asymptotics of the distributions of $N$ and $A$.
Thus, for example,
the error terms in Theorem \ref{strongclt} could have instead been
stated in terms of the variances.

For completeness, we should mention what is known about the growth of
(say) $\mathbb E[A]$,
which can be effectively estimated using elementary geometric and
combinatorial techniques. In dimension two, one has
%
%
\begin{equation}\label{planegrowth}
\log[\Area(K)]\ll\mathbb E[A]\ll[\Area(K)]^{1/3}.
\end{equation}
[In particular, the error terms in Theorem \ref{strongclt} go to zero
as $\Area(K)\to\infty$.]
The estimate (\ref{planegrowth}) is a consequence of the \textit{economic
cap covering lemma} of B\'ar\'any
and Larman \cite{baranycap} in combination with other estimates in
\cite
{baranycap} and those of Groemer \cite{groemersphere} (in
fact, their results apply to higher dimensions as well). We remark that
the lower asymptotic is achieved when $K$ is a
polygon, and the upper asymptotic is achieved when $\partial K$ is
$C^2$ with nonvanishing curvature.

We conclude by remarking that in recent years there has been
significant progress in the study of random \textit{polytopes}, but again
most results deal only with the cases when
(i) $K$ is a polytope \cite{baranypolytope}, and (ii) $\partial K$ is
$C^2$ with nonvanishing Gauss curvature \cite{reitzner1,vu1}. We
believe that an approach similar to ours should be
possible in higher dimensions as well. This would shed new light on
problems in that setting, and ultimately show that there is no qualitative
difference between the cases (i) and (ii).

\subsection{The uniform model random polygons}

A model related to $\Pi_K$ is $P_{K,n}:=\mbox{conv. hull.}(X_1,\ldots,X_n)$ where $X_i$ are i.i.d. uniformly
in $K$; $P_{K,n}$ is called a \textit{random polygon}. This is often
referred to as the ``uniform model'' whereas $\Pi_K$ is
the ``Poisson model.'' Morally they are the same process in the limit
$\Area(K)=n\to\infty$ (though making this precise is
often difficult). It has been a well-known open
problem to prove central limit theorems for functionals of $P_{K,n}$.
For instance, Van Vu \cite{vuopen} has asked the question
of whether a central limit theorem holds for $A(P_{K,n})$, though the
problem is a very natural one in the study of random polygons,
a subject that began with work of R\'enyi and Sulanke \cite
{renyisulanke1,renyisulanke2}. Theorems \ref{strongclt} and \ref{strongvarest}
both carry over to the setting of $P_{K,n}$, thus answering this
question in the affirmative.
\begin{corollary}\label{uniformclt}
As $n\to\infty$, we have the following central limit theorems for $P_{K,n}$:
%
%
\begin{eqnarray}
\sup_x\biggl|P\biggl(\frac{N-\mathbb E[N]}{\sqrt{\Var N}}\leq
x
\biggr)-\Phi(x)\biggr|&\to&0,
\\
%
%
\sup_x\biggl|P\biggl(\frac{A-\mathbb E[A]}{\sqrt{\Var A}}\leq
x
\biggr)-\Phi(x)\biggr|&\to&0
\end{eqnarray}
uniformly over all convex $K$. Here $\Phi(x)=P(Z\leq x)$ where $Z$ is
the standard normal distribution.
\end{corollary}
\begin{corollary}\label{uniformvarest}
As $n\to\infty$, we have the following estimates for $P_{K,n}$:
%
%
\begin{equation}
\mathbb E[N]\asymp\Var N\asymp\frac n{\Area(K)}\mathbb E[A]\asymp
\biggl(\frac n{\Area(K)}\biggr)^2\Var A
\end{equation}
uniformly over all convex $K$.
\end{corollary}

As in the case of the Poisson model, these results are well known in
the field in the two cases (i) $K$ a polygon and
(ii) $\partial K$ of class $C^2$ with nonvanishing curvature. The
innovation in this paper is that all $K$ are
treated uniformly.

A detailed derivation of Corollaries \ref{uniformclt} and \ref
{uniformvarest} from Theorems \ref{strongclt} and \ref{strongvarest}
will appear elsewhere \cite{pardonsupplement}.
Suffice it to say here that they are almost immediate consequences of
the corresponding results on the Poisson model once
one proves that when $n=\Area(K)$, the variables $N(P_{K,n})$ and
$N(\Pi
_K)$ [as well as $A(P_{K,n})$ and $A(\Pi_K)$]
have the same expectation and variance up to a small enough error.

\section{The basic decomposition}\label{sec2}

In this section, we illustrate our basic approach. We will aim for
Theorem \ref{strongclt}, and Theorem \ref{strongvarest}
will be a corollary of our methods.

\begin{figure}

\includegraphics{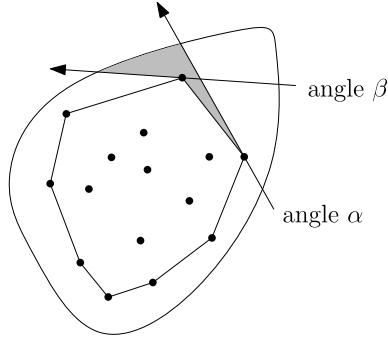}

\caption{Illustration of $A(\alpha,\beta)$.}\label{Adecomp}
\end{figure}

First, we observe that the functionals $N$ and $A$ both enjoy
decompositions into local pieces. We define $N(\alpha,\beta)$ to equal
the number of edges of $\Pi$ whose angle lies in the interval $[\alpha
,\beta]\subset\mathbb R/2\pi$. The definition of $A(\alpha,\beta)$
is best explained graphically (see Figure \ref{Adecomp}). Thus for any
fixed sequence of angles $\alpha_1<\alpha_2<\cdots<\alpha_L$, we have
the following decompositions:
%
%
\begin{eqnarray}
N&=&N(\alpha_1,\alpha_2)+\cdots+N(\alpha_L,\alpha_1),
\\
%
%
A&=&A(\alpha_1,\alpha_2)+\cdots+A(\alpha_L,\alpha_1).
\end{eqnarray}
During the proof, we often do not need to distinguish between whether
we are dealing with $N$ or $A$. Thus we will use $X(\Pi)$
to denote either $N$ or $A$ when a statement holds for both.

%
%

A central limit theorem will follow if we can find a choice of $\{
\alpha
_i\}$ such that the moments of $X(\alpha_i,\alpha_{i+1})$ are bounded
uniformly,
and such that the dependence between $X(\alpha_i,\alpha_{i+1})$ and
$X(\alpha_j,\alpha_{j+1})$ becomes small as $|i-j|\to\infty$. Our
construction
is to choose $\{\alpha_i\}$ so that the intervals $[\alpha_i,\alpha
_{i+1}]$ have constant \textit{affine invariant measure} (a measure
depending on $K$).
In this paper, we give a more or less explicit description of the
affine invariant measure, which in practice should allow its easy
estimation for any given class
of convex sets, and thus a complete description of the behavior of
random Poisson polygons and random polygons. As we remarked in the
\hyperref[sec1]{Introduction}, a key
result is the compactness of various configuration spaces.

After fixing notation in Section \ref{notdefsec}, we define the affine
invariant measure in Section \ref{affinvmdefsec}. Section \ref{compactsec}
is devoted to the crucial step of proving the compactness of the
configuration spaces. Using the information coming from compactness:
\begin{itemize}
\item In Section \ref{momentsec}, we estimate the moments of $X$
(Proposition \ref{finitemomentgen}).
\item In Section \ref{independsec}, we estimate the long range
dependence of $X$ (Proposition \ref{exponentialdecreaseaffinvmeasure}).
\item In Section \ref{varianceestimatesec}, we recall an estimate the
variance of $X$ due to Imre B\'ar\'any and Matthias Reitzner
(Proposition \ref{weakvarest}).
\end{itemize}
The remainder of the paper contains the explicit deduction of Theorems
\ref{strongclt} and~\ref{strongvarest}.

\section{Notation and definitions}\label{notdefsec}

In this paper, $K$ will always denote a (boun\-ded) convex set in
$\mathbb R^2$.

We warn the reader that in most of the literature, one fixes $\Area
(K)=1$ and then considers
a Poisson process of intensity $\lambda\to\infty$. We have chosen
instead to use the normalization
$\lambda=1$ and let $\Area(K)\to\infty$. This is convenient for us
because it makes many of our formulas simpler to state.

Any constants implied by the symbols $\ll$, $\gg$ or $\asymp$ are
absolute; in particular they are not allowed to depend on $K$.
There will be times when we require $\Area(K)\gg1$; this is no real
restriction to us since in the end we will take $\Area(K)\to\infty$.
The group $\Aff(2)=\mathbb R^2\rtimes\SL_2(\mathbb R)$ is the group of
(oriented) area preserving affine transformations of $\mathbb R^2$;
it acts naturally on the entire problem studied here.

%
\begin{figure}

\includegraphics{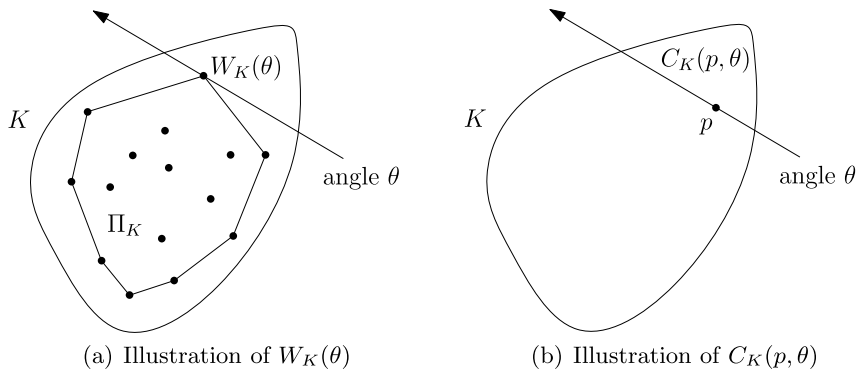}

\caption{Illustration of some definitions.}\label{basic}
\end{figure}

Many of the following definitions are illustrated in Figure \ref
{basic}. We may leave out the subscript $K$ later when doing so is unambiguous.
\begin{definition}
We define the random variable $W_K(\theta)$ to be the vertex of $\Pi_K$
which has an oriented tangent line at angle $\theta$.
This is illustrated in Figure \ref{basic}(a).
\end{definition}
\begin{definition}
A \textit{cap at angle $\theta$} is the intersection of $K$ with a
half-plane $H_\theta$ at angle $\theta$. We
may specify a cap at angle $\theta$ by giving either its area $r$ or a
point $p\in\partial H_\theta$.
These are denoted $C_K(r,\theta)$ and $C_K(p,\theta)$, respectively; the
latter is illustrated in Figure \ref{basic}(b).
\end{definition}
\begin{definition}
We define the real number $A_K(p,\theta)$ to be the area of the cap
$C_K(p,\theta)$.
\end{definition}
\begin{lemma}\label{Wdistr}
The random variable $W_K(\theta)$ has probability distribution given by
$\exp(-A_K(p,\theta)) \,dp$ where $dp$ is the Lebesgue measure.
\end{lemma}
\begin{pf}
This follows directly from the definition of a Poisson point process.
\end{pf}
\begin{definition}
We define the function $f_K(x,\theta)\dvtx[0,1]\times\mathbb R/2\pi\to
\mathbb R$ as follows:
%
%
\begin{equation}
f_K(x,\theta)=\cases{\mbox{length of $(\partial H_\theta)\cap K$,}\cr
\hspace*{10.38pt}\qquad\mbox{where $\displaystyle C_K\biggl(\log\frac1x,\theta\biggr)=H_\theta\cap K$,}
\vspace*{2pt}\cr
\hspace*{32.75pt}\mbox{if }x>\exp(-\Area(K)),\cr
0,\qquad \mbox{if $x\leq\exp(-\Area(K))$}.}
\end{equation}
\end{definition}

It will be important to have the following bound on the growth of $f$:
\begin{lemma}\label{boundonF}
If $y\leq x$, then
%
%
\begin{equation}
\frac{f(y)}{\sqrt{-\log y}}\leq\frac{f(x)}{\sqrt{-\log x}}.
\end{equation}
\end{lemma}

The bound above is sharp; for instance $f(x)=\operatorname{const}\cdot
\sqrt{-\log x}$ for $K=\{x,y\geq0\}$ (i.e., the first quadrant).

\begin{pf*}{Proof of Lemma \ref{boundonF}}
Project $K$ along the lines at angle $\theta$ to get a height function
$h\dvtx[0,\infty)\to\mathbb R_{\geq0}$;
in Figure \ref{figureH}, $h(\ell)$ is the length of the thick segment.
%
%
\begin{figure}

\includegraphics{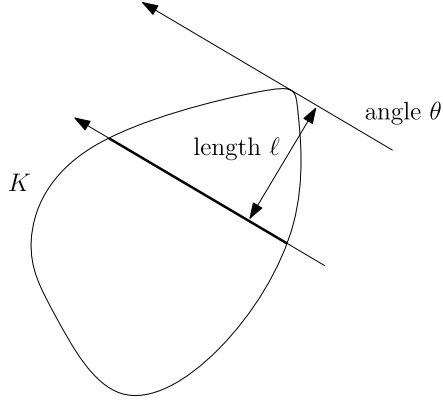}

\caption{Illustration of the function $h$.}\label{figureH}
\end{figure}
Now if
$A(\ell)=\int_0^\ell h(\ell') \,d\ell'$ then $f(\exp(-A(\ell
)))=h(\ell
)$. Thus
we see that it suffices to show that the function
%
%
\begin{equation}
\frac{h(\ell)}{\sqrt{A(\ell)}}
\end{equation}
is decreasing. Differentiating with respect to $\ell$, we see that it
suffices to show that
%
%
\begin{equation}
h(\ell)^2-2h'(\ell)A(\ell)\geq0.
\end{equation}
For $\ell=0$, the left-hand side is clearly nonnegative, and the
derivative of the
left-hand side equals $-2h''(\ell)A(\ell)$, which is $\geq0$ by
concavity of $h$.
\end{pf*}
\begin{lemma}\label{fnotwild}
If $\Area(K)\geq2\log\frac1x$, then $f(y)\leq2 f(x)$ for $y\geq x$.
\end{lemma}
\begin{pf}
Refer to Figure \ref{figureCompr}. The area of the upper trapezoid is
$\leq\log\frac1x$
%
%
\begin{figure}

\includegraphics{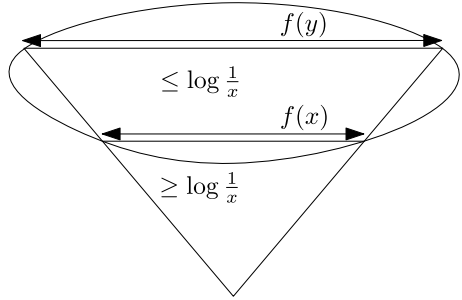}

\caption{Illustration of an inequality.}\label{figureCompr}
\end{figure}
since it is contained in $C(\log\frac1x,\theta)$. The area of the
lower triangle is $\geq\log\frac1x$
since it contains $K\setminus C(\log\frac1x,\theta)$ and $\Area
(K)\geq
2\log\frac1x$. Similar triangles gives the following inequality:
%
%
\begin{equation}
\frac{f(y)-f(x)}{\log(1/x)}\leq\frac{f(x)}{\log(1/x)}.
\end{equation}
Simplifying yields $f(y)\leq2f(x)$.
\end{pf}

\section{The affine invariant measure}\label{affinvmdefsec}

\begin{proposition}
For every $g\in\Aff(2)$, we have
%
%
\begin{equation}
r_g^\ast[f_{gK}(x,\theta)^2 \,d\theta]=f_{K}(x,\theta
)^2\,d\theta,
\end{equation}
where $r_g\dvtx\mathbb R/2\pi\to\mathbb R/2\pi$ is the action of $g$ on
line slopes. We say \textup{``}$f(x,\theta
)^2\,d\theta$ is affine invariant.\textup{''}
\end{proposition}
\begin{pf}
Define $\mathbf v(\theta)$ to be the vector of length $f(x,\theta)$
parallel to the chord whose length gives $f(x,\theta)$.
Then we have
%
%
\begin{equation}
\int_{\theta_1}^{\theta_2}f(x,\theta)^2 \,d\theta=\int_{\theta
_1}^{\theta
_2}\mathbf v(\theta)\times d\mathbf v(\theta).
\end{equation}
The right-hand side is invariant under the action of $\Aff(2)$, so the
result follows.
\end{pf}
\begin{definition}\label{affineinvariantmeasuredef}
We define the \textit{affine invariant measure} to be $\mu
_K:=f_K(e^{-1},\theta)^2 \,d\theta$.
\end{definition}

The $\varepsilon$-\textit{wet part} of $K$ is defined as the union of all
caps of area $\varepsilon$. In the literature,
estimates for random polygons are frequently expressed in terms of the
area of the $\varepsilon$-wet part of $K$.
It is, perhaps, not surprising that our notion of the affine invariant
measure is related to the area of the wet part in
the following manner:
\begin{lemma}\label{wetpartaffinemeasure}
One has the following relation:
%
%
\begin{equation}
\Area\biggl(\bigcup_{\gamma\in[\alpha,\beta]}C_K(1,\gamma)
\biggr)=1+\frac18\mu_K([\alpha,\beta]).
\end{equation}
\end{lemma}
\begin{pf}
Consider the area swept out by the line segments bounding the caps of
area $1$ at angles $\gamma\in[\alpha,\beta]$ (area covered
twice is counted twice). On the one hand, this area just equals
%
%
\begin{equation}\label{areaeasy}
2\Area\biggl(\bigcup_{\gamma\in[\alpha,\beta]}C_K(1,\gamma
)\biggr)-\Area
(C_K(1,\alpha))-\Area(C_K(1,\beta)).
\end{equation}
On the other hand, we may express the area as an integral $d\theta$.
Each line segment rotates about its midpoint (since
the area of the caps is constant), so the area covered is just the
$d\theta$ integral of
$\int_{-f(e^{-1},\theta)/2}^{f(e^{-1},\theta)/2}|y| \,dy=\frac
14f(e^{-1},\theta)^2$. Comparing this with (\ref{areaeasy}) yields
the result.
\end{pf}

\section{Compactness of configuration spaces}\label{compactsec}

\begin{definition}
Define a configuration space $\mathcal C(r)$ for $r>0$ as follows. The
objects of $\mathcal C(r)$
are convex subsets of $\mathbb R^2$ of area $r$ with a distinguished
line segment on their boundary. As a set,
$\mathcal C(r)$ is equal to everything of the form $(H\cap K,(\partial
H)\cap K)$, where $K$ is any convex
set of area $\geq2r$ and $H$ is a half-plane such that $H\cap K$ has
area $r$. A typical member of $\mathcal C(r)$ is illustrated in Figure
\ref{capfigure}(a).
We emphasize that the space $\mathcal C(r)$ does not depend on any
choice of convex set $K$; rather it
is the space of \textit{all} caps of area $r$ that come from some convex
set of area $\geq2r$.

%
%
\begin{figure}[b]

\includegraphics{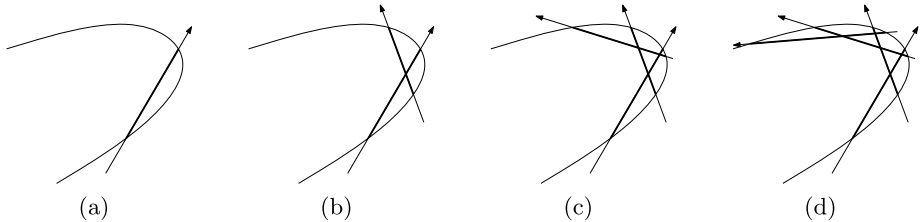}

\caption{A series of caps.}\label{capfigure}
\end{figure}

We call $\mathcal C(r)$ the configuration space of caps of area $r$. If
$c\in\mathcal C(r)$,
then we call the distinguished part of its boundary its \textit{flat
boundary} and the undistinguished part of its boundary its \textit{convex
boundary}.
We let the \textit{half-plane of} $c$ equal the unique half-plane which
contains $c$ and whose boundary contains the flat boundary of $c$
(this is exactly the $H$ appearing above).
\end{definition}

We topologize $\mathcal C(r)$ by using the Hausdorff metric to compare
both the set and its distinguished
subset. Explicitly, $d((A,A_0),(B,B_0))=d(A,B)+d(A_0,B_0)$. Let us
observe that there is a natural action
of $\Aff(2)$ on $\mathcal C(r)$; it is continuous. Certainly $\mathcal
C(r)$ is not compact, since the group
$\Aff(2)$ is noncompact. However, we will show directly that $\mathcal
C(r)/\Aff(2)$ is compact. This simple
fact will be an essential tool in virtually all of the estimates in the
remainder of this paper.
\begin{lemma}\label{singlecompact}
The space $\mathcal C(r)/\Aff(2)$ is compact.
\end{lemma}
\begin{pf}
Let $c_1,c_2,\ldots$ be a sequence of elements of $\mathcal C(r)/\Aff
(2)$. Pick
representatives $\tilde c_1,\tilde c_2,\ldots$ in $\mathcal C(r)$ so
that the flat
part of $\partial\tilde c_i$ is the unit line segment on the $x$-axis,
$\tilde c_i$ is contained in
the upper half-plane, and the highest $y$-coordinate of any point in
$\tilde c_i$ is attained at $(\frac12,h_i)$.
This is illustrated in Figure \ref{figureCompact}.

%
%
\begin{figure}

\includegraphics{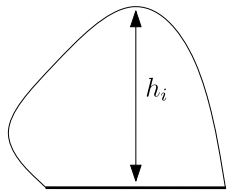}

\caption{Compactness of $\mathcal C(r)/\Aff(2)$.}\label{figureCompact}
\end{figure}

By Lemma \ref{fnotwild}, we conclude that $\Area(K)\geq2r$ implies
that every
horizontal chord in $\tilde c_i$ has length $\leq2$. This implies that
$-\frac32\leq x\leq\frac52$ for any $x$-coordinate
of a point in $\tilde c_i$. On the other hand, $\tilde c_i$ contains a
triangle of base $1$ and height~$h_i$, so by comparing
areas we must have $\frac12h_i\leq r$. Thus we conclude that $\tilde
c_i\subseteq[-\frac32,\frac52]\times[0,2r]$. It is well known that
the space of
convex sets of fixed volume in some bounded region of $\mathbb R^d$
given the Hausdorff topology is compact (this is the so-called Blaschke
selection theorem). Thus
we conclude that there exists a subsequence of $\tilde c_i$ that converges.
\end{pf}
\begin{definition}
We define the complex configuration space $\mathcal C(r_1,\varepsilon
,r_2)$ for
$r_1,r_2>0$ and $0<\varepsilon<\min(r_1,r_2)$ as follows. We let
$\mathcal
C(r_1,\varepsilon,r_2)$ denote
a particular subset of $\mathcal C(r_1)\times\mathcal C(r_2)$. An
ordered pair $(c_1,c_2)\in\mathcal C(r_1)\times\mathcal C(r_2)$
is in $\mathcal C(r_1,\varepsilon,r_2)$ if and only if it satisfies the
following:
\begin{itemize}
\item$\Area(c_1\cap c_2)=\varepsilon$.
\item If $H_1$ is the half-plane of $c_1$, then $H_1\cap c_2=c_1\cap c_2$.
\item If $H_2$ is the half-plane of $c_2$, then $c_1\cap H_2=c_1\cap c_2$.
\item It holds that $\operatorname{angle}(H_1)<\operatorname
{angle}(H_2)<\operatorname
{angle}(H_1)+\pi$.
\end{itemize}
We then give $\mathcal C(r_1,\varepsilon,r_2)$ the subspace topology.

One can see that the middle two conditions taken together just mean
that $c_1$ and $c_2$ coincide on $H_1\cap H_2$,
and the last condition just says that $c_1$ precedes $c_2$ if we
traverse their convex boundary counterclockwise.
Examples appear in Figure \ref{capfigure}(b) and in Figure \ref{figureMoment}.
\end{definition}
\begin{lemma}\label{doublecompact}
The space $\mathcal C(r_1,\varepsilon,r_2)/\Aff(2)$ is compact.
\end{lemma}
\begin{pf}
Let $(c_1,d_1),(c_2,d_2),\ldots$ be a sequence of elements of the
quotient $\mathcal C(r_1,\varepsilon,r_2)/\Aff(2)$. Lift
these to a sequence $(\tilde c_1,\tilde d_1),(\tilde c_2,\tilde
d_2),\ldots$ in $\mathcal C(r_1,\varepsilon,r_2)$ where we
assume (after passing to a subsequence using Lemma \ref{singlecompact})
that $\tilde c_1,\tilde c_2,\ldots$
is convergent to $\tilde c\in\mathcal C(r_1)$.

%
%
\begin{figure}

\includegraphics{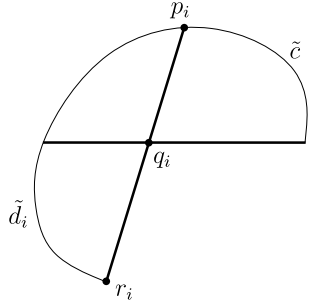}

\caption{Compactness of $\mathcal C(r_1,\varepsilon,r_2)/\Aff
(2)$.}\label{figureCompactD}
\end{figure}

Now refer to Figure \ref{figureCompactD}. Label the intersection of the
flat boundary of $\tilde d_i$
with the convex boundary of $\tilde c_i$ as $p_i$. Label the
intersection of the flat boundaries of
$\tilde d_i$ and $\tilde c_i$ as $q_i$. Label the intersection of the
flat boundary of $\tilde d_i$
with its convex boundary other than $p_i$ as $r_i$. Clearly we can
extract a subsequence for which $p_i$
converges to a point $p$ on the convex boundary of $\tilde c$, and then
extract a further subsequence for which $q_i$
converges to a point $q$ on the flat boundary of $\tilde c$. The only
subtlety in this proof is to observe
that $0<\varepsilon<r_1$ shows that $p$ and $q$ are not on the corners of
$\tilde c$.

Given $p$ and $q$, the boundedness of the area of $\tilde d_i$
implies that $r_i$ is bounded, so we extract another subsequence for
which additionally $r_i$ converges to a point $r$.
Now it is easy to see that the fixing of $\tilde c,p,q,r$ provide only
a bounded set for $\tilde d_i$ to range over,
so compactness follows again using the Blaschke selection theorem.
\end{pf}
\begin{lemma}\label{findsequence}
There exists an absolute constant $M_0<\infty$ such that if we are
given $K$ and angles $\alpha<\beta$ with $\mu_K([\alpha,\beta
])\geq
M_0$, then we can find
a sequence $\alpha\leq\gamma_0<\gamma_1<\cdots<\gamma_L\leq\beta
$ so
that $(C_K(\gamma_{i-1},1),C_K(\gamma_i,1))\in\mathcal C(1,\frac12,1)$
and $L\asymp\mu_K([\alpha,\beta])$.
\end{lemma}
\begin{pf}
Let $\gamma_0=\alpha$. Now define $\gamma_i$ inductively for $i\geq1$
as follows. The function
%
%
\begin{equation}
\Area\bigl(C(1,\gamma_{i-1})\cap C(1,\gamma)\bigr) \qquad\mbox{for $\gamma\in
[\gamma
_{i-1},\gamma_{i-1}+\pi]$}
\end{equation}
is strictly decreasing until it reaches zero, where it remains
constant. Thus there exists a unique $\gamma_i$ so
that $\Area(C(1,\gamma_{i-1})\cap C(1,\gamma_i))=\frac12$. We now have
an infinite chain of angles
$\alpha=\gamma_0<\gamma_1<\gamma_2<\cdots$ so that
$C(1,\gamma_i)\cap C(1,\gamma_{i+1})$ has area $\frac12$ for $i\geq
0$. This is illustrated in Figure \ref{capfigure}.

Let $L$ be the maximum index such that $\gamma_L\leq\beta$. Note that
since $\mathcal C(1,\frac12,1)$ is
compact, there exist absolute constants $0<Y_1<Y_2<\infty$ (not
depending on $K$) such that
%
%
\begin{equation}
Y_1<\mu_K([\gamma_i,\gamma_{i+1}])<Y_2
\end{equation}
for all $i$. Thus we conclude that
%
%
\begin{equation}
Y_1L<\mu_K([\gamma_0,\gamma_L])\leq\mu_K([\alpha,\beta])<Y_2(L+1),
\end{equation}
which is sufficient.
\end{pf}

\section{A moment estimate}\label{momentsec}

An ingredient in the central limit theorems for the polygonal case is a
moment estimate \cite{groen1}, page 341, Lemma 2.5, and \cite{groen2},
page 36, Lemma 2.1. Here, we
prove an analogous estimate in general.
\begin{proposition}\label{finitemomentgen}
Let $X$ denote either $N$ or $A$. There exist absolute constants
$M_0<\infty$ and $\varepsilon>0$ such that for any convex $K$ and interval
$[\alpha,\beta]$
with $\mu_K([\alpha,\beta])\geq M_0$, we have the following estimate:
%
%
\begin{equation}
\mathbb E\exp(\lambda X_K(\alpha,\beta))\ll1 \qquad\mbox{for all
}|\lambda|<\varepsilon/\mu_K([\alpha,\beta]).
\end{equation}
\end{proposition}
\begin{pf}
We can split up $[\alpha,\beta]$ into subintervals of small affine
invariant measure, and use Cauchy's inequality,
%
%
\begin{equation}
\mathbb E \exp(\lambda[A+B])\leq\sqrt{[\mathbb E \exp(2\lambda
A)][\mathbb E \exp(2\lambda B)]},
\end{equation}
so it suffices to show that
there exist $\delta>0$ and $\varepsilon>0$ so that for all $K$ and
$[\alpha
,\beta]$ satisfying $\mu_K([\alpha,\beta])\leq\delta$, it holds that
the moment generating function $\mathbb E \exp(\lambda X_K(\alpha
,\beta
))$ is $\ll1$ for all $|\lambda|<\varepsilon$.

Since $\mathcal C(1,\frac12,1)$ is compact, the affine invariant
measure of the interval between the angles
of $c_1$ and $c_2$ is bounded below. Thus we conclude that it suffices
to show that for every $(c_1,c_2)\in\mathcal C(1,\frac12,1)$,
the moment generating function of $X_K(\alpha,\beta)$ is defined in a
neighborhood of zero where $\alpha$ is the angle of
$c_1$ and $\beta$ is the angle of~$c_2$.

Now we may put such an element $(c_1,c_2)\in\mathcal C(1,\frac12,1)$
in a standard position in $\mathbb R^2$
by requiring that both boundary segments have equal length, and that
the angles of $c_1$ and $c_2$ are
$0$ and $\frac\pi2$, respectively, (see Figure \ref{figureMoment}).

%
%
\begin{figure}

\includegraphics{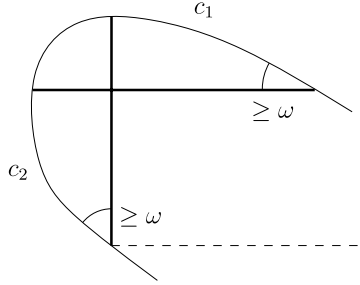}

\caption{Two adjacent caps.}\label{figureMoment}
\end{figure}

Thus, given the configuration in Figure \ref{figureMoment}, we would
like to show that for sufficiently
small $\lambda>0$, we have $\mathbb E \exp(\lambda X_K(0,\frac\pi
2))\ll1$. First, write
%
%
\begin{eqnarray}
&&\mathbb E\exp\biggl(\lambda X_K\biggl(0,\frac\pi2\biggr)\biggr)\nonumber\\[-8pt]\\[-8pt]
&&\qquad=\int_K\mathbb E\biggl[\exp
\biggl(\lambda
X_K\biggl(0,\frac\pi2\biggr)\biggr)\Big|W(0)=p\biggr]
\,dP\bigl(W(0)=p\bigr).\nonumber
\end{eqnarray}

If $X=N$, then $X_K(0,\frac\pi2)$ is bounded by the number of points of the
Poisson process in the region $C(W(0),\frac\pi2)\setminus C(W(0),0)$.
An elementary calculation shows that $\mathbb E \exp(\lambda\Xi
(k))=\exp
(k[e^\lambda-1])$, where $\Xi(k)$ is a Poisson distribution of
parameter~$k$.
We may assume $|\lambda|<1$, so $e^\lambda-1<2|\lambda|$. Thus in
this case
%
%
\begin{eqnarray}
&&\mathbb E\biggl[\exp\biggl(\lambda X_K\biggl(0,\frac\pi2\biggr)\biggr)\Big|W(0)=p\biggr]\nonumber\\[-8pt]\\[-8pt]
&&\qquad\leq\exp\biggl(2|\lambda
|\Area\biggl(C\biggl(p,\frac\pi2\biggr)\Bigm\backslash
C(p,0)\biggr)\biggr).\nonumber
\end{eqnarray}

If $X=A$, then $X_K(0,\frac\pi2)$ is bounded by $C(W(0),\frac\pi
2)\setminus C(W(0),0)$, so we have
%
%
\begin{eqnarray}
&&\mathbb E\biggl[\exp\biggl(\lambda X_K\biggl(0,\frac\pi2\biggr)\biggr)\Big|W(0)=p\biggr]\nonumber\\[-8pt]\\[-8pt]
&&\qquad\leq\exp\biggl(|\lambda
|\Area
\biggl(C\biggl(p,\frac\pi2\biggr)\Bigm\backslash
C(p,0)\biggr)\biggr).\nonumber
\end{eqnarray}

Thus in both cases, we have the estimate
%
%
\begin{eqnarray}\quad
&&\mathbb E\exp\biggl(\lambda X_K\biggl(0,\frac\pi2\biggr)\biggr)\nonumber\\[-8pt]\\[-8pt]
&&\qquad\leq\int_K\exp\biggl(2|\lambda
|\Area
\biggl(C\biggl(p,\frac\pi2\biggr)\Bigm\backslash C(p,0)\biggr)\biggr)\exp(-A(p,0))
\,dp\nonumber
\end{eqnarray}
recalling Lemma \ref{Wdistr}.

By compactness of $\mathcal C(1,\frac12,1)/\Aff(2)$, the angle where
the convex part of $c_i$ meets the flat boundary of $c_i$
is bounded below by an absolute constant (say by~$\omega$, see Figure
\ref{figureMoment}). Similarly, the lengths of the flat
parts of $c_1$ and $c_2$ are bounded above absolutely (say by $R\geq
1$). Thus the area above the dotted line in Figure \ref{figureMoment}
is bounded above absolutely, say by $B=2+R^2+R^2\cot\omega$.

Now we claim that
%
%
\begin{equation}
\Area\biggl(C\biggl(p,\frac\pi2\biggr)\Bigm\backslash C(p,0)\biggr)\leq B+f(p,0)^2\cot\omega
\end{equation}
[recall that $f(p,0)$ is the length of $\ell\cap K$ where $\ell$ is the
horizontal line passing through $p$].
If $p\in c_1$, then the area of $C(p,\frac\pi2)\setminus C(p,0)$ is
$\leq B$ by definition. If $p\notin c_1$,
then argue as follows: the area of $C(p,\frac\pi2)\setminus C(p,0)$
above the dotted line is certainly less
than $B$, and the area of $C(p,\frac\pi2)\setminus C(p,0)$ below the
dotted line is bounded by $f(p,0)^2\cot\omega$.

Thus we have
%
%
\begin{eqnarray}
&&\mathbb E \exp\biggl(\lambda X_K\biggl(0,\frac\pi2\biggr)\biggr)\nonumber\\[-8pt]\\[-8pt]
&&\qquad\leq e^{2|\lambda|B}\int
_K\exp
(2|\lambda|f(p,0)^2\cot\omega)\exp(-A(p,0)) \,dp.\nonumber
\end{eqnarray}
If we substitute $x=\exp(-A(p,0))$, then the integral becomes
%
%
\begin{equation}
\mathbb E \exp\biggl(\lambda X_K\biggl(0,\frac\pi2\biggr)\biggr)\leq e^{2|\lambda|B}\int
_0^1\exp(2|\lambda|f(x,0)^2\cot\omega) \,dx.
\end{equation}
Now $f(e^{-1},0)\leq R$, so $f(x,0)\leq R\sqrt{-\log x}$ for $x\leq
e^{-1}$ by Lemma \ref{boundonF}, and
$\Area(K)\gg1$ implies $f(x,0)\leq2R$ for $x\geq e^{-1}$ by Lemma
\ref{fnotwild}. Thus we conclude that
%
%
\begin{eqnarray}
\mathbb E \exp\biggl(\lambda X_K\biggl(0,\frac\pi2\biggr)\biggr)&\leq& e^{2|\lambda|B}\int
_0^{e^{-1}}x^{-2|\lambda|R^2\cot\omega} \,dx\nonumber\\[-8pt]\\[-8pt]
&&{}+e^{2|\lambda|B}\int
_{e^{-1}}^1e^{8|\lambda|R^2\cot\omega} \,dx,\nonumber
\end{eqnarray}
which is bounded absolutely for small enough $|\lambda|$.
\end{pf}

\section{A dependence estimate}\label{independsec}

\begin{definition}
If $S\subset\mathbb R/2\pi$ is an interval, then we let $\mathcal
F_S^{(K)}$ be the $\sigma$-algebra which keeps track of $W_K(\theta)$
for $\theta\in S$.
\end{definition}

For example, $\Pi_K$ is $\mathcal F_S^{(K)}$-measurable if and only if
$S=\mathbb R/2\pi$.

The type of dependence estimate we prove will be an $\alpha$-mixing
estimate, that is, an estimate on
$|P(A\cap B)-P(A)P(B)|$ where $A$ and $B$ are events that are supposed
to be almost independent. This
type of estimate has been used previously in studying random polygons;
we were motivated to prove our
estimate by a similar result in \cite{groen1}, page 341, Theorem 2.3.
\begin{lemma}\label{independencethm}
Let $[\theta_1,\theta_2]$ and $[\psi_1,\psi_2]$ be two disjoint
intervals in $\mathbb R/2\pi$.
Let $A\in\mathcal F_{[\theta_1,\theta_2]}$ and $B\in\mathcal
F_{[\psi
_1,\psi_2]}$. Then
%
%
\begin{eqnarray}\label{simplealpha}
&&|P(A\cap B)-P(A)P(B)|\nonumber\\[-8pt]\\[-8pt]
&&\qquad\ll\sum_{i,j\in\{1,2\}}\int_K\exp
(-A(p,\theta_i))\exp(-A(p,\psi_j)) \,dp.\nonumber
\end{eqnarray}
\end{lemma}

The proof is an elementary calculation and is given in the
\hyperref[app]{Appendix}. The object of this section is to reexpress
the right-hand side of (\ref{simplealpha}) in terms of the affine
invariant measure.
\begin{lemma}\label{helpindp}
There exists an absolute constant $\delta>0$ such that if $\theta\leq
\psi\leq\theta+\pi$, then area of
$C(1,\theta)\cap C(1,\psi)$ is $\ll\exp(-\delta\mu_K([\theta
,\psi]))$.
\end{lemma}

\begin{pf}
We use Lemma \ref{findsequence} to construct a sequence $\theta
=\gamma
_0<\gamma_1<\cdots<\gamma_L\leq\psi$ so that
$\Area(C(\gamma_i,1)\cap C(\gamma_{i+1},1))=\frac12$ and $L\asymp
\mu
_K([\theta,\psi])$. From this decomposition,
we see that it suffices to show that there exists $\delta>0$ such that
for all $i$
%
%
\begin{equation}
\Area\bigl(C(\gamma_i,1)\cap C(\gamma_0,1)\bigr)\leq(1-\delta)\Area
\bigl(C(\gamma
_{i-1},1)\cap C(\gamma_0,1)\bigr).
\end{equation}
Now we know that $C(\gamma_0,1)=K\cap H$ for some half-plane $H$ and
that additionally
$\Area(C(\gamma_{i-1},1)\cap C(\gamma_0,1))=\Area(C(\gamma
_{i-1},1)\cap
H)\leq\frac12$.
Hence it suffices to show that
%
%
\begin{equation}
\Area\bigl(C(\gamma_i,1)\cap H\bigr)\leq(1-\delta)\Area\bigl(C(\gamma
_{i-1},1)\cap H\bigr),
\end{equation}
whenever $\Area(C(\gamma_{i-1},1)\cap H)\leq\frac12$ and
$\operatorname
{angle}(H)\in(\gamma_i-\pi,\gamma_{i-1})$.

%
%
\begin{figure}[b]

\includegraphics{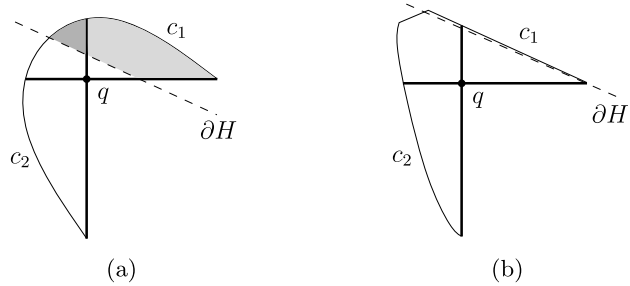}

\caption{Intersecting caps.}\label{figureExp}
\end{figure}

Remember that $C(\gamma_i,1)$ and $C(\gamma_{i-1},1)$ have intersection
$\frac12$. Thus it suffices to show that for
every $(c_1,c_2)\in\mathcal C(1,\frac12,1)$, the following is true:
%
%
\begin{equation}\label{expdecrease}
\frac{\Area(c_2\cap H)}{\Area(c_1\cap H)}<1-\delta,
\end{equation}
whenever $\Area(c_1\cap H)\leq\frac12$ and $\operatorname
{angle}(H)\in
(\operatorname{angle}(c_2)-\pi,\operatorname{angle}(c_1))$
[see Figure \ref{figureExp}(a)]. Here, if we put $c_1$ and $c_2$
in standard position (i.e., as in Figure \ref{figureExp}, with both
flat boundaries of equal length), then $\partial H$
has negative slope. Denote by\vspace*{1pt} $q$ the intersection of the flat
boundaries of $c_1$ and $c_2$. Then since $\Area(c_1\cap H)\leq\frac12$,
we must have $q\notin H$. From this, we see that $c_2\cap H\subseteq
c_1\cap H$, so we may rewrite (\ref{expdecrease}) as
%
%
\begin{equation}\label{expdecrease2}
\frac{\Area((c_1\setminus c_2)\cap H)}{\Area(c_2\cap H)}>\delta.
\end{equation}
The minimum of this expression is clearly a continuous function on
$\mathcal C(1,\frac12,1)$, and is by definition invariant
under the action of $\Aff(2)$. We know that $\mathcal C(1,\frac
12,1)/\Aff(2)$ is compact, so it suffices to show that for any fixed
configuration $(c_1,c_2)$, expression (\ref{expdecrease2}) is bounded
below away from zero. Certainly, if this ratio were
approaching zero, then $\Area((c_1\setminus c_2)\cap H)\to0$. However
in this case, the situation is illustrated in Figure
\ref{figureExp}(b), where it is clear that ratio (\ref{expdecrease2})
in fact does \textit{not} approach zero, but rather
some appropriate ratio of lengths of the boundaries of the caps. Thus
we are done.
\end{pf}
\begin{lemma}
There exists an absolute constant $\delta>0$ such that if $\theta\leq
\psi\leq\theta+\pi$
%
%
\begin{equation}
\int_K\exp(-A(p,\theta))\exp(-A(p,\psi)) \,dp\ll\exp(-\delta\mu
_K([\theta
,\psi])).
\end{equation}
\end{lemma}
\begin{pf}
We pick the unique $\theta_1,\psi_1$ so
that $\theta<\theta_1<\psi_1<\psi$ and $\mu_K([\theta$, $\theta
_1])=\mu
_K([\theta_1,\psi_1])=\mu_K([\psi_1,\psi])$.

Define
%
%
\begin{equation}
S_p=C(p,\theta)\cup C(p,\psi)=\bigcup_{\theta\leq\alpha\leq\psi
}C(p,\alpha),
\end{equation}
so $\Area(S_p)\leq A(p,\theta)+A(p,\psi)$.

Now if $A(p,\alpha)\geq1$ for
all $\alpha\in[\theta,\theta_1]$, then by Lemma \ref{wetpartaffinemeasure},
the area of $S_p$ is $\gg\mu_K([\theta,\theta_1])=\frac13\mu
_K([\theta
,\psi])$. The same applies
if $A(p,\alpha)\geq1$ for $\alpha\in[\psi_1,\psi]$. Thus in both of
these cases, we
conclude that $A(p,\theta)\gg\mu_K([\theta,\psi])$ or $A(p,\psi
)\gg\mu
_K([\theta,\psi])$.

If $A(p,\theta_2)<1$ for some $\theta_2\in[\theta,\theta_1]$ and
$A(p,\psi_2)<1$ for some $\psi_2\in[\psi_1,\psi]$, then
necessarily $p\in C(1,\theta_1)\cap C(1,\psi_1)$. Thus we know that for
all $p\in K$, at least one of
the following is true:
\begin{itemize}
\item$p\in C(1,\theta_1)\cap C(1,\psi_1)$,
\item$A(p,\theta)\gg\mu_K([\theta,\psi])$,
\item$A(p,\psi)\gg\mu_K([\theta,\psi])$.
\end{itemize}
By elementary integration, the integral over the second and third
regions is $\ll\exp(-\delta\mu_K(\theta,\psi))$.
The area of the first region is $\ll\exp(-\delta\mu_K(\theta,\psi
))$ by Lem\-ma~\ref{helpindp}, so we are done.
\end{pf}
\begin{proposition}\label{exponentialdecreaseaffinvmeasure}
There exists an absolute constant $\delta>0$ so that if $[\theta
_1,\theta_2]$ and $[\psi_1,\psi_2]$ are two
disjoint intervals in $\mathbb R/2\pi$, and we have events $A\in
\mathcal F_{[\theta_1,\theta_2]}$ and $B\in\mathcal F_{[\psi_1,\psi
_2]}$, then
%
%
\begin{equation}
|P(A\cap B)-P(A)P(B)|\ll\sum_{i,j\in\{1,2\}}\exp
(-\delta
\mathfrak d_K(\theta_i,\psi_i)),
\end{equation}
where $\mathfrak d_K(\alpha,\beta)$ denotes $\mu_K([\alpha,\beta
])$ if
$\alpha\leq\beta\leq\alpha+\pi$ and
$\mu_K([\beta,\alpha])$ if instead $\beta\leq\alpha\leq\beta
+\pi$.
\end{proposition}

The reader may wonder exactly what follows from an $\alpha$-mixing
estimate. We won't answer that here, though
we will record here two lemmas that will be useful later whose
hypotheses are $\alpha$-mixing estimates.
\begin{lemma}[(\cite{centrallimittheorembook}, page 115, Lemma 1(6))]\label
{alphahassmallcov}
Suppose $X$ and $Y$ are random variables taking values in $\mathbb R$
such that
%

\begin{equation}
|P(X\in A\ \&\ Y\in B)-P(X\in A)P(Y\in B)|<\alpha
\end{equation}
for all $A,B\subseteq\mathbb R$. Then we have
%
%
\begin{equation}
|{\Cov}(X,Y)|\leq6(\mathbb E|X|^3)^{1/3}(\mathbb E|Y|^3)^{1/3}\alpha^{1/3}.
\end{equation}
\end{lemma}
\begin{lemma}\label{alphagood}
Suppose $X$ and $Y$ are random variables taking values in $\mathbb R$
such that
%
%
\begin{equation}
|P(X\in A\ \&\ Y\in B)-P(X\in A)P(Y\in B)|<\alpha
\end{equation}
for all $A,B\subseteq\mathbb R$. Let $Z=X+Y$, and let $\tilde Z$ equal
the sum of independent copies of $X$ and $Y$. Then we have
%
%
\begin{equation}
{\sup_x}|P(Z\leq x)-P(\tilde Z\leq x)|\ll\sqrt\alpha.
\end{equation}
\end{lemma}
\begin{pf}
Let $-\infty=x_0<x_1<\cdots<x_N=\infty$ be any finite increasing
sequence of real numbers. Then we have
%
%
\begin{eqnarray}
P(Z\leq0)&\geq&\sum_{i=1}^NP\bigl(X\in(x_{i-1},x_i]\ \&\ Y\leq-x_i\bigr)
\nonumber\\[-8pt]\\[-8pt]
&\geq&-N\alpha+\sum_{i=1}^NP\bigl(X\in(x_{i-1},x_i]\bigr)P(Y\leq
-x_i).\nonumber
\end{eqnarray}
Now using the definition of $\tilde Z$, we can bound this below by
%
%
\begin{eqnarray}
P(Z\leq0)&\geq&-N\alpha+P(\tilde Z\leq0)\nonumber\\[-8pt]\\[-8pt]
&&{}-\sum_{i=1}^NP\bigl(X\in
(x_{i-1},x_i)\bigr)P\bigl(Y\in
(-x_i,-x_{i-1})\bigr).\nonumber
\end{eqnarray}
Thus we find that
%
%
\begin{eqnarray}\label{alphaprobs}
&&
P(Z\leq0)-P(\tilde Z\leq0)\nonumber\\[-8pt]\\[-8pt]
&&\qquad\geq-N\alpha-\sum_{i=1}^NP\bigl(X\in
(x_{i-1},x_i)\bigr)P\bigl(Y\in(-x_i,-x_{i-1})\bigr).\nonumber
\end{eqnarray}
Now choose $K-1$ real numbers $-\infty=u_0<u_1<\cdots<u_K=\infty$ so
that the probability that $X$ falls in the open interval $(u_{i-1},u_i)$
is $\leq K^{-1}$ for all $i$. Do the same for $Y$ to get $v_i$'s. Then
let the $x_i$'s be the union of the $u_i$'s and $-v_i$'s
(so $N\leq2K$). With this choice, we see that each of the
probabilities in the last sum of (\ref{alphaprobs}) is $\leq K^{-1}$, so
their product is $\leq K^{-2}$. Hence the right-hand side is $\geq
-2K\alpha-2KK^{-2}$.
Now choosing $K$ to equal the nearest integer to $\alpha^{-1/2}$, we
conclude that $P(Z\leq0)-P(\tilde Z\leq0)\geq-\mbox{const}\cdot
\alpha^{1/2}$.
By a symmetric\vspace*{1pt} argument, we get the other inequality, so $|P(Z\leq
0)-P(\tilde Z\leq0)|\ll\alpha^{1/2}$, which is
sufficient.
\end{pf}

\section{A variance estimate}\label{varianceestimatesec}

The task of providing a lower bound on the variance of $N$ and $A$ has
already been completed by B\'ar\'any and
Reitzner~\cite{baranypreprint}, page 4, Theorem 2.1.
They prove the following theorem.
\begin{proposition}\label{weakvarest}
Provided $\mu_K([\alpha,\beta])\gg1$, we have the estimates
%
%
\begin{eqnarray}
\Var N(\alpha,\beta)&\gg&\mu_K([\alpha,\beta]),
\\
%
%
\Var A(\alpha,\beta)&\gg&\mu_K([\alpha,\beta]).
\end{eqnarray}
\end{proposition}

In fact, B\'ar\'any and Reitzner's result is valid for random polytopes
as well. They only state this estimate
in the case $[\alpha,\beta]=\mathbb R/2\pi$, though their proof is
valid in general. We also
note that they phrase their result in terms of the area of the
$\varepsilon
$-wet part of $K$; we have replaced this with the
affine invariant measure using Lem\-ma~\ref{wetpartaffinemeasure}.

\section{\texorpdfstring{Proof of Theorem \protect\ref{strongvarest}}{Proof of Theorem 1.2.}}

Let $X$ denote either $N$ or $A$.

From linearity of the expectation, one immediately observes that
$\mathbb E[X]\asymp\mu_K(\mathbb R/2\pi)$.
Proposition \ref{weakvarest} implies that
%
%
\begin{equation}
\mu_K(\mathbb R/2\pi)\ll\Var X.
\end{equation}
Thus it suffices to show the reverse inequality. For this, simply
decompose $\mathbb R/2\pi$ into $L$ intervals of affine
invariant measure $\asymp1$, and then write
%
%
\begin{eqnarray}
\Var X&=&\sum_{i=1}^L\Var X(\alpha_i,\alpha_{i+1})\nonumber\\[-8pt]\\[-8pt]
&&{}+2\sum_{1\leq
i<j\leq
L}\Cov(X(\alpha_i,\alpha_{i+1}),X(\alpha_j,\alpha_{j+1})).\nonumber
\end{eqnarray}
Proposition \ref{finitemomentgen} shows that the sum of variances is
$\ll L$.
Proposition \ref{exponentialdecreaseaffinvmeasure},
Lem\-ma~\ref{alphahassmallcov} and Proposition \ref{finitemomentgen}
imply that the sum of covariances is $\ll L$. Hence the right-hand side
is $\ll L$ as needed.

\section{\texorpdfstring{Proof of Theorem \protect\ref{strongclt}}{Proof of Theorem 1.1.}}\label{proofsection}

We need the following central limit theorem appearing in a survey
article by Sunklodas \cite{sunklodas}:
\begin{theorem}[(In English translation \cite{centrallimittheorembook}, page 133, Theorem 10)]\label{generalCLT}
Let $X=\sum_{i=1}^LX_i$ where $X_1,\ldots,X_L$ are random variables.
Additionally suppose that:
\begin{itemize}
\item$\mathbb E |X_i|^3\leq C_1$,
\item$X_1,\ldots,X_L$ are $\alpha$-mixing with $\alpha\leq C_2\exp
(-\delta|i-j|)$,
\end{itemize}
for some $\delta>0$ and $C_1,C_2<\infty$. Then there exists $M<\infty$
such that
%
%
\begin{equation}
\sup_{x\in\mathbb R}\biggl|P\biggl(\frac{X-\mathbb E X}{\sqrt{\Var
X}}\leq x\biggr)-\Phi(x)\biggr|\leq M\frac{L(\log L)^2}{(\Var X)^{-3/2}}.
\end{equation}
\end{theorem}

We have everything necessary to apply Theorem \ref{generalCLT} to $N$
and $A$, except that our decomposition
is ``circular.'' Thus, for example, Theorem \ref{generalCLT} shows
immediately that $N(\alpha,\alpha+\pi)$ satisfies a
central limit theorem for any $\alpha$, but does not directly apply to
give a central limit theorem for $N$.
For completeness, we include the following proof, where we derive
Theorem \ref{strongclt} just using Theorem \ref{generalCLT}
as a black box. The reader who is willing to believe the natural
extension of Theorem~\ref{generalCLT} to our situation may
want to omit it, as it is essentially just a straightforward calculation.

\begin{pf*}{Proof of Theorem \ref{strongclt}}
Suppose $K$ is given with $\Area(K)\gg1$. Let $X$ denote either $N$ or
$A$. In this proof $\delta>0$ denotes some
positive absolute constant, possibly different at each occurrence.

The function $f(\alpha)=\mu_K([\alpha,\alpha+\pi])$ on $\mathbb
R/2\pi$
satisfies $f(\alpha)+f(\alpha+\pi)=\mu_K(\mathbb R/2\pi)$.
Thus by continuity we may find $\alpha$ such that $\mu_K([\alpha
,\alpha
+\pi])=\mu_K([\alpha+\pi,\alpha+2\pi])$. Without loss of generality,
we may assume $\mu_K([0,\pi])=\mu_K([\pi,2\pi])$. Set $L=\mu
_K(\mathbb
R/2\pi)$.

We let $\ell$ denote a quantity much smaller than $L$ (we will
eventually let $\ell$ equal some large multiple of $\log L$).
We pick $\alpha_1,\beta_1,\alpha_2,\beta_2$ so that $\mu([0,\alpha
_1])=\mu([\beta_1,\pi])=\mu([\pi,\alpha_2])=\mu([\beta_2,2\pi
])=\ell$.
Then we set
%
%
\begin{eqnarray}
X_1&=&X(\alpha_1,\beta_1),\\
X_2&=&X(\alpha_2,\beta_2).
\end{eqnarray}
Observe that by partitioning $[\alpha_i,\beta_i]$ into intervals of
affine invariant measure $\asymp1$, we may apply Theorem
\ref{generalCLT} (appealing to Propositions \ref{finitemomentgen} and
\ref{exponentialdecreaseaffinvmeasure}). Thus remembering
Proposition \ref{weakvarest}, we may write
%
%
\begin{equation}\label{xiclose}
\biggl|P\biggl(\frac{X_i-\mathbb E[X_i]}{\sqrt{\Var X_i}}\leq x
\biggr)-\Phi(x)\biggr|\ll\frac{\log^2L}{\sqrt L}.
\end{equation}
Let $\tilde Y$ equal the sum of independent copies of $X_1$ and $X_2$.
Then (\ref{xiclose}) implies that
%
%
\begin{equation}
\biggl|P\biggl(\frac{\tilde Y-\mathbb E[\tilde Y]}{\sqrt{\Var\tilde
Y}}\leq x\biggr)-\Phi(x)\biggr|\ll\frac{\log^2L}{\sqrt L}.
\end{equation}

By Proposition \ref{exponentialdecreaseaffinvmeasure}, $X_1$ and $X_2$
are $\alpha$-mixing with $\alpha\ll e^{-\delta\ell}$.
Proposition \ref{finitemomentgen} shows $\mathbb E[|X_i|^3]^{1/3}\ll
L$. If we let $Y=X_1+X_2$, then Lemma \ref{alphagood}
implies that $|P(Y\leq x)-P(\tilde Y\leq x)|\ll e^{-\delta\ell}$. Lemma
\ref{alphahassmallcov} implies
%
%
\begin{equation}
\Cov(X_1,X_2)\ll L^2e^{-\delta\ell}.
\end{equation}
Since $\Var Y=\Var\tilde Y+2\Cov(X_1,X_2)$ and $\Var\tilde Y\asymp L$,
we have $\Var Y=(1+O(Le^{-\delta L}))\Var\tilde Y$.
Hence $|\Phi(\sqrt{\Var Y}x)-\Phi(\sqrt{\Var\tilde Y}x)|\ll
Le^{-\delta
\ell}$. Hence we conclude that
%
%
\begin{eqnarray}\label{Yisgood}
\biggl|P\biggl(\frac{Y-\mathbb E[Y]}{\sqrt{\Var Y}}\leq x
\biggr)-\Phi
(x)\biggr|&\ll&\frac{\log^2L}{\sqrt L}+e^{-\delta\ell}+Le^{-\delta
\ell}
\nonumber\\[-8pt]\\[-8pt]
&\ll&\frac{\log^2L}{\sqrt L}+Le^{-\delta\ell}.
\nonumber
\end{eqnarray}

Now the final part of our argument is to translate this into a
statement about $X$. Let $E=X(\beta_2,\alpha_1)+X(\beta_1,\alpha_2)$.
Thus by definition, we have $X=Y+E$. Using Proposition \ref
{finitemomentgen}, it is evident that $\mathbb E[\exp(\delta\ell
^{-1}E)]\ll1$ for
some absolute $\delta>0$. From this, we conclude that $P(\exp(\delta
\ell
^{-1}E)\geq M)\ll M^{-1}$. Thus $P(E\geq M)\ll e^{-\delta M/\ell}$.
Now we pick $M=\ell^2$, so that
%
%
\begin{equation}\label{Egood}
P(E\geq\ell^2)\ll e^{-\delta\ell}.
\end{equation}
Now examine (\ref{Yisgood}), and consider what this says about $P
(\frac{Y+E-\mathbb E[Y]}{\sqrt{\Var Y}}\leq x)$. We have
$\sqrt{\Var Y}\asymp\sqrt L$, so (\ref{Egood}) implies that
$P(E/\sqrt
{\Var Y}\notin[0,\ell^2/\sqrt L])\ll e^{-\delta\ell}$. Thus
%
%
\begin{eqnarray}
&&\biggl|P\biggl(\frac{Y+E-\mathbb E[Y]}{\sqrt{\Var Y}}\leq x
\biggr)-\Phi
(x)\biggr|
\nonumber\\
&&\qquad\ll\biggl|P\biggl(\frac{Y-\mathbb E[Y]}{\sqrt{\Var Y}}\leq x
\biggr)-\Phi
(x)\biggr|+e^{-\delta\ell}+\frac{\ell^2}{\sqrt L}
\\
&&\qquad\ll\frac{\log^2L}{\sqrt L}+Le^{-\delta\ell}+\frac{\ell^2}{\sqrt
L}.\nonumber
\end{eqnarray}
Now $\mathbb E[E]\asymp\ell$, so adding $\mathbb E[E]$ in the numerator
adds at most $\ell/\sqrt L$ to the error. Hence
%
%
\begin{equation}
\biggl|P\biggl(\frac{X-\mathbb E[X]}{\sqrt{\Var Y}}\leq x
\biggr)-\Phi
(x)\biggr|\ll\frac{\log^2L}{\sqrt L}+Le^{-\delta\ell}+\frac{\ell
^2}{\sqrt L}.
\end{equation}
Observe that $\Var E\ll\ell$, so $\Var X=\Var Y+2\Cov(Y,E)+\Var
E=\Var
Y+O(\sqrt L\sqrt\ell)+O(\ell)$, so the relative
error is $\ll\frac{\sqrt\ell}{\sqrt L}+\frac\ell L$. Thus we have
%
%
\begin{eqnarray}
&&\biggl|P\biggl(\frac{X-\mathbb E[X]}{\sqrt{\Var X}}\leq x
\biggr)-\Phi
(x)\biggr|\nonumber\\[-8pt]\\[-8pt]
&&\qquad\ll\frac{\log^2L}{\sqrt L}+Le^{-\delta\ell}+\frac{\ell
^2}{\sqrt L}+\frac{\sqrt\ell}{\sqrt L}+\frac\ell L.\nonumber
\end{eqnarray}
Taking $\ell$ to equal a sufficiently large multiple of $\log L$, we
achieve the desired estimate.
\end{pf*}

\vspace*{-18pt}

\begin{appendix}\label{app}

\section*{\texorpdfstring{Appendix: Proof of Lemma
\protect\lowercase{\ref{independencethm}}}{Appendix: Proof of Lemma 7.2.}}

Lemmas \ref{lem1} and \ref{lem2} below combine easily to give Lemma
\ref{independencethm}. The proof of Lemma \ref{lem1} follows \cite{groen1},
where similar manipulations are performed.
\setcounter{theorem}{0}
\begin{lemma}\label{lem1}
Let $[\theta_1,\theta_2]$ and $[\psi_1,\psi_2]$ be two disjoint
intervals in $\mathbb R/2\pi$.
Let $A\in\mathcal F_{[\theta_1,\theta_2]}$ and $B\in\mathcal
F_{[\psi
_1,\psi_2]}$. Then
%
%
\setcounter{equation}{0}
\begin{eqnarray}
&&|P(A\cap B)-P(A)P(B)|
\nonumber\\[-8pt]\\[-8pt]
&&\qquad\leq2\sum_{i,j\in\{1,2\}}\iint_{(p,q)\in R(\theta_i,\psi
_j)}dP\bigl(W(\theta_i)=p\bigr) \,dP\bigl(W(\psi_j)=q\bigr),
\nonumber
\end{eqnarray}
where $R(\alpha,\beta)$ is the set of pairs $(p,q)\in K\times K$ such
that it is impossible that $W(\alpha)=p$ and $W(\beta)=q$.
\end{lemma}
\begin{pf}
We have that $P(A\cap B)$ is given by
%
%
\begin{eqnarray}\qquad
&&\iiiint_{K^4}P\bigl(A\cap B|(W(\theta_1),W(\theta_2),W(\psi_1),W(\psi
_2))=(p_1,p_2,q_1,q_2)\bigr)\nonumber\\[-8pt]\\[-8pt]
&&\qquad\hspace*{18.25pt}{}\times dP\bigl((W(\theta_1),W(\theta_2),W(\psi_1),W(\psi_2))
=(p_1,p_2,q_1,q_2)\bigr)
\nonumber
\end{eqnarray}
and $P(A)P(B)$ by
%
%
\begin{eqnarray}\qquad
&&\iint_{K^2}P\bigl(A|(W(\theta_1),W(\theta_2))=(p_1,p_2)\bigr) \,dP\bigl((W(\theta
_1),W(\theta_2))=(p_1,p_2)\bigr)
\nonumber\\
&&\qquad{}\times\iint_{K^2}P\bigl(B|(W(\psi_1),W(\psi_2))=(q_1,q_2)\bigr)\\
&&\hspace*{59pt}{}\times dP\bigl((W(\psi_1),W(\psi_2))=(q_1,q_2)\bigr).
\nonumber
\end{eqnarray}
Now given $W(\theta_1)$, $W(\theta_2)$, $W(\psi_1)$ and $W(\psi_2)$,
the events $A$ and $B$
are independent. In other words the two integrands above are equal
(although the measures are not). Hence we conclude that
%
%
\begin{eqnarray}\label{tobound}
\hspace*{8pt}&&|P(A\cap B)-P(A)P(B)|\hspace*{-15pt}\nonumber\\
\hspace*{8pt}&&\qquad\leq\iiiint_{K^4}
\bigl|dP\bigl((W(\theta_1),W(\theta_2),W(\psi_1),W(\psi_2))=(p_1,p_2,q_1,q_2)\bigr)\\
\hspace*{8pt}&&\hspace*{4pt}\quad\qquad{}-dP\bigl((W(\theta_1),W(\theta_2))=(p_1,p_2)\bigr) \,dP\bigl((W(\psi_1),W(\psi
_2))=(q_1,q_2)\bigr)\bigr|.\nonumber
\end{eqnarray}
Now define the set
%
%
\begin{eqnarray}
&&
R_{\theta_1,\theta_2,\psi_1,\psi_2}
=\{(p_1,p_2,q_1,q_2)\in K^4\dvtx\nonumber\\
&&\hspace*{69.2pt}(W(\theta_1),W(\theta_2),W(\psi_1),W(\psi
_2))=(p_1,p_2,q_1,q_2)\\
&&\hspace*{232.5pt}\mbox{ is impossible}\}.
\nonumber
\end{eqnarray}
We will calculate the right-hand side of (\ref{tobound}) by splitting
up the integral as $I(R_{\theta_1,\theta_2,\psi_1,\psi
_2})+I(R_{\theta
_1,\theta_2,\psi_1,\psi_2}^\complement)$
(i.e., the\vspace*{2pt} integral over $R_{\theta_1,\theta_2,\psi_1,\psi_2}$ and
the integral over its complement). Since the first measure in question
$dP((W(\theta_1),W(\theta_2)$, $W(\psi_1),W(\psi_2))=(p_1,p_2,q_1,q_2))$
is supported on $R_{\theta_1,\theta_2,\psi_1,\psi_2}^\complement$,
we trivially have that
%
%
\begin{eqnarray}\label{hiddeneq}
\hspace*{4pt}
&&\iiiint_{R_{\theta_1,\theta_2,\psi_1,\psi_2}^\complement}
\bigl[dP\bigl((W(\theta_1),W(\theta_2),W(\psi_1),W(\psi
_2))=(p_1,p_2,q_1,q_2)\bigr)\nonumber\\
\hspace*{4pt}
&&\hspace*{42pt}{}-dP\bigl((W(\theta_1),W(\theta_2))=(p_1,p_2)\bigr) \,dP\bigl((W(\psi_1),W(\psi
_2))=(q_1,q_2)\bigr)\bigr]
\nonumber\\[-8pt]\\[-8pt]
\hspace*{4pt}
&&\qquad=\iiiint_{R_{\theta_1,\theta_2,\psi_1,\psi_2}}
dP\bigl((W(\theta_1),W(\theta_2))=(p_1,p_2)\bigr) \nonumber\\
\hspace*{4pt}
&&\qquad\hspace*{83pt}{}\times dP\bigl((W(\psi_1),W(\psi
_2))=(q_1,q_2)\bigr).\nonumber
\end{eqnarray}
Now observe that on $R_{\theta_1,\theta_2,\psi_1,\psi
_2}^\complement
$, we have
%
%
\begin{eqnarray}
&&dP\bigl((W(\theta_1),W(\theta_2),W(\psi_1),W(\psi
_2))=(p_1,p_2,q_1,q_2)\bigr)\nonumber\\[-8pt]\\[-8pt]
&&\qquad\geq
dP\bigl((W(\theta_1),W(\theta_2))=(p_1,p_2)\bigr) \,dP\bigl((W(\psi_1),W(\psi
_2))=(q_1,q_2)\bigr).\hspace*{-31pt}
\nonumber
\end{eqnarray}
From this, it is clear that equation (\ref{hiddeneq}) is equivalent to
%
%
\begin{equation}
I(R_{\theta_1,\theta_2,\psi_1,\psi_2}^\complement)=I(R_{\theta
_1,\theta_2,\psi_1,\psi_2}).
\end{equation}
Thus the right-hand side of (\ref{tobound}) in fact equals
$2I(R_{\theta
_1,\theta_2,\psi_1,\psi_2})$. Hence
we conclude that $|P(A\cap B)-P(A)P(B)|$ is bounded
above by
%
%
\begin{eqnarray}
&&2\iiiint_{R_{\theta_1,\theta_2,\psi_1,\psi_2}}
dP\bigl((W(\theta_1),W(\theta_2))=(p_1,p_2)\bigr)\nonumber\\[-8pt]\\[-8pt]
&&\qquad\hspace*{56.9pt}{}\times dP\bigl((W(\psi_1),W(\psi
_2))=(q_1,q_2)\bigr).\nonumber
\end{eqnarray}
If we define $R(\theta,\psi):=\{(p,q)\in K^2\dvtx(W(\theta),W(\psi))=(p,q)$
is impossible$\}$, then
%
%
\begin{equation}
R_{\theta_1,\theta_2,\psi_1,\psi_2}=\bigcup_{i,j\in\{1,2\}
}R(\theta
_i,\psi_j)\times K^2.
\end{equation}
Thus we conclude that
%
%
\begin{eqnarray}\qquad
&&|P(A\cap B)-P(A)P(B)|\nonumber\\
&&\qquad\leq2\sum_{i,j\in\{1,2\}}\iiiint
_{R(\theta_i,\psi_j)\times K^2}
dP\bigl((W(\theta_1),W(\theta_2))=(p_1,p_2)\bigr)\\
&&\qquad\hspace*{130pt}{}\times dP\bigl((W(\psi_1),W(\psi
_2))=(q_1,q_2)\bigr).
\nonumber
\end{eqnarray}
Integrating out the undesired indices on the right-hand side yields the
correct result.
\end{pf}
\begin{lemma}\label{lem2}
We have
%
%
\begin{eqnarray}
&&\iint_{R_{\alpha,\beta}}dP\bigl(W(\alpha)=p\bigr) \,dP\bigl(W(\beta)=q\bigr)
\nonumber\\[-8pt]\\[-8pt]
&&\qquad\leq2\int_K\exp(-A(p,\alpha))\exp(-A(p,\beta)) \,dp.
\nonumber
\end{eqnarray}
\end{lemma}
\begin{pf}
This relies on the observation that $R_{\alpha,\beta}=\{(p,q)\in
K^2\dvtx q\notin H_{p,\beta}\}\cup\{(p,q)\in K^2\dvtx p\notin
H_{q,\alpha}\}$.
Now recalling that $dP(W(\alpha)=p)=A(p,\alpha) \,dp$ (Lemma \ref
{Wdistr}), we calculate
%
%
\begin{eqnarray}
&&\int_K\biggl[\int_{K-H_{p,\beta}}dP\bigl(W(\beta)=q\bigr)\biggr]\,dP\bigl(W(\alpha)=p\bigr)
\nonumber\\
&&\qquad=\int_K\exp(-A(p,\beta)) \,dP\bigl(W(\alpha)=p\bigr)\\
&&\qquad=\int_K\exp(-A(p,\beta))\exp(-A(p,\alpha)) \,dp.\nonumber
\end{eqnarray}
Thus the result follows.
\end{pf}
\end{appendix}

\section*{Acknowledgments}

The author thanks Yakov Sinai for conversations about this work as well
as opportunities to speak about it in seminars. Thanks also go
to Imre B\'ar\'any and Matthias Reitzner for posting \cite
{baranypreprint} on B\'ar\'any's website. The author is also grateful
to the referee
for many helpful comments.

\printaddresses

\end{document}